\documentclass[a4page, 11pt]{article}

\usepackage{fullpage}
\usepackage{graphics}
\usepackage{pdfpages}
\usepackage{indentfirst}
\usepackage{setspace}
\usepackage{natbib}
\usepackage{verbatim}
\usepackage{lineno}
\usepackage[font=small]{caption}
\modulolinenumbers[5]

\newcommand{\compresslist}{%
	\setlength{\itemsep}{1.5pt}%
	\setlength{\parskip}{0.2pt}%
	\setlength{\parsep}{0.2pt}%
}

\begin{document}

\title{Contrasting Prediction Methods for Early Warning Systems at Undergraduate Level}
\author{Emma Howard, Maria Meehan \& Andrew Parnell}
\maketitle

\begin{abstract}
In this study, we investigate prediction methods for an early warning system for a large STEM undergraduate course. Recent studies have provided evidence in favour of adopting early warning systems as a means of identifying at-risk students. Many of these early warning systems rely on data from students' engagement with Learning Management Systems (LMSs). Our study examines eight prediction methods, and investigates the optimal time in a course to apply an early warning system. We present findings from a statistics university course which has a large proportion of resources on the LMS Blackboard and weekly continuous assessment. We identify weeks 5-6 of our course (half way through the semester) as an optimal time to implement an early warning system, as it allows time for the students to make changes to their study patterns whilst retaining reasonable prediction accuracy. Using detailed (fine-grained) variables, clustering and our final prediction method of BART (Bayesian Additive Regressive Trees) we are able to predict students' final grade by week 6 based on mean absolute error (MAE) to 6.5 percentage points. We provide our R code for implementation of the prediction methods used in a GitHub repository.

\end{abstract}

\singlespacing

\section{Introduction} 

Early warning systems to identify at-risk students (of dropping out or failing) are in practical use in large classes and online courses \citep*{Corrigan2015, Wolff2014, Pistilli2012}. We provide findings from a large first year statistics course in which most of the learning materials are available online and therefore student engagement with them can be measured via the LMS Blackboard. We acknowledge the impact course design, in particular weekly continuous assessment, has on developing early warning systems. We contrast results from eight prediction methods (Random Forest; BART; XGBoost; Principal Components Regression; Support Vector Machine; Neural Network; Multivariate Adaptive Regression Splines; and K-Nearest Neighbours) and the impact of cluster membership (based on student engagement) on reducing prediction error. We reasonably predict a student's final grade as early as week five of a twelve week semester. This study was completed using R software and we have provided our R code on GitHub at https://github.com/ehoward1/Early-Warning-System, and in the appendix.

This study forms part of a larger goal to use the predictions we create to allow for more precisely targeted interventions for poorly performing students. Determining the timing at which these interventions should occur is one of the key goals of this study. We would like to intervene as early as possible, but with little information from the LMS and necessarily limited continuous assessment at the start of the semester the predictions are inaccurate. This accuracy increases as we move through the semester but at the price of intervening later and so lessening the impact of any interventions. We monitor the performance of the predictive models on a week by week cumulative basis. For each week, we aim to predict the final percentage mark of the student based on all current information. We do not dichotomise students’ performance to pass/fail unlike many other studies \citep*{Marbouti2016, Azcona2015, Macfadyen2010} which would lessen the accuracy. At week 6 (of a 12 week semester) we obtain a mean absolute error (MAE) of approximately 6.5 percentage points.

The structure of our paper is as follows: in Section 2 we discuss the rationale and prediction methods behind current early warning systems. In Section 3 we outline our approach to developing an accurate prediction method for an early warning system. We extend current research on the development of early warning systems through: using `new' prediction methods including BART; identifying an `optimal time'; and including cluster membership. In Section 4, we discuss the data analytics decisions made and present the results for our course \emph{Practical Statistics}. Finally we progress to the discussion and conclusion of these results in Section 5.

\section{Previous Work on Early Warning Systems}
\subsection{Prediction Modelling for Early Warning Systems}

In this section, we examine the stages in creating a prediction model for an early warning system (detailed data collection; variable selection; prediction modelling; and clustering). Advancements in learning interfaces allow for fine-grained collection of data. \citet{Azcona2015} highlight that ``fine-grained (microscopic) analytics data should yield better results than coarse-grained (macroscopic)'' (p. 589). An example of a coarse-grained variable is total count of resources accessed online. In comparison fine-grained data analytics refers to extracting each log entry of a student, and all the information it contains for example: the number of slides visited; number of successful compilations; and time spent on platform \citep{Azcona2015}. Their argument is that through using more detailed variables, more powerful prediction models can be created. However this trades off against simplicity; simple models with a small number of variables are easier to interpret and understand.

Variables based on students' demographic/historic data, continuous assessment results and LMS usage have been collected for early warning systems \citep*{You2016, Pistilli2012}. LMS data can include length of time on a LMS system, number of visits to a module page, contributions to a module discussion thread et cetera. Depending on the prediction models selected, the dataset is reduced to a small number of `important' variables. 

There are numerous types of prediction models used for learning analytics. \citet*{Ga2016} note that researchers have produced prediction models by using classification algorithms such as EM, C4.5, Naive Bayes Classifier, and Support Vector Machines. Logistic regression and multiple regression modelling are often used as prediction models \citep{Macfadyen2010, Waddington2016}, with logistic regression being considered the most popular prediction method for educational settings \citep{Marbouti2016}. Hierarchical mixed models \citep*{Joksimovic2015, You2016}, K-nearest neighbor \citep{Marbouti2016}, neural network models \citep*{Calvo-Flores2006}, and decision tree methods \citep{Azcona2015} are also methods employed.  A common use of prediction models in learning analytics is to identify whether a student will pass or fail the course based on the binary response variable `pass/fail'. The use of a binary response variable dichotomises students' performance percentage marks. Studies using binary response variable include \citet{Azcona2015, Macfadyen2010, Calvo-Flores2006}.

A key point to note is that predictive models are usually applied to a single course rather than used for several courses. \cite*{Pantucek2013} propose that this may be because each course is structured differently, and therefore dictates what learners are doing. \citet{Ga2016} investigate generalised predictive models that can be applied to multiple courses, however they note that the inherent differences in disciplines cause specific variables to be strong for some courses, and weak for other courses. Hence, the nature of the course should be considered before selecting variables for an early warning system. \citet{Ga2016} believe ``the understanding of practical needs in specific instructional and learning contexts is the primary driver for the development and deployment of learning analytics methods'' (p. 83).  

Clustering also plays a significant role in learning analytics through its ability to identify students' engagement levels or learning strategies statistically. When investigating a blended course \citet*{Lust2011} identify three patterns of tool-use using k-means clustering: the no-users; the intensive users; and the incoherent users. In comparison, \citet*{Howard2017} used model-based clustering to identify four clusters of behavioural engagement in a large mathematics course where students have the choice to use lectures or/and online videos. \cite*{White2016} use Latent Class Analysis to identify four clusters of engagement in a large blended business course. In their discussion they identify what resources each cluster engaged with, and when during the semester these resources were engaged with. 

\subsection{Early Warning Systems in Practice}

One of the best known examples of an early warning system is in Purdue University \citep*{Sclater2016, Ferguson2012, Pistilli2012} who introduced `Course Signals' (CS) or a `traffic light system' whereby students can see whether they are likely to succeed in their course based on a traffic light colour on their learner interface. For example a green colour indicates a high likelihood of succeeding. This prediction of success is based on prediction models using all available student background information and LMS interactions. If a student is identified as at-risk, the lecturer has the option of providing corrective measures including: posting of a traffic signal indicator on the student's CMS home page; sending e-mail messages or reminders; sending text messages; referring the student to an academic advisor or academic resource centre; or organising a face-to-face meeting. \cite{Pistilli2010} found that the results of their interventions (based on a control group versus an experimental group) were: students seeking help earlier; lower D's and F's recorded; more B's and C's; and students felt more than a `number', that is less isolated. Other benefits of Course Signals discussed by \cite{Sclater2016} are students using the subject help desks more, and greater attendance at additional tutorials.  

One prime reason for the implementation of an early warning system is to detect students at-risk of dropping out of courses. \cite{Pistilli2010} state that most early warning systems rely on midterm grades reported by lecturers. By the time midterms have been corrected it is often far into the semester, and students may have already dropped out. It is crucial that early warning systems operate in the early stages of the semester. However, a balance has to be achieved with the accuracy of the model. As the methods, models, variables and response variable used in identifying at-risk students vary from study to study, it remains difficult to contrast the studies and identify which study has obtained the most accurate results. Results are impacted by the truncating of students' performances to the binary pass/fail variable. Dichotomizing is usually performed for simplicity however this can lead to: lower accuracy through loss of valuable information; a decrease in the predictive power; and in general there is a risk of getting results that may not make sense \citep*{Fedorov2009, Royston2006}. Many studies have reported results of identifying at-risk students at the end of the course/semester however for early warning systems this is impractical. Ideally we wish to support all students from the beginning of the semester. For a prediction model, the beginning of the semester is too early to identify at-risk students. For early warning systems, a balance needs to be obtained between the increasing accuracy of the system and the diminishing impact of intervening as we move through the semester. In this paper we refer to the balance between the two as the `optimal time'.  

\subsection{Research Questions}

Our study aims to explore developing a prediction model for an early warning system taking into account the benefits of cluster analysis. Furthermore our study aims to identify an `optimal time' in the semester when an early warning system could be implemented. Hence our research questions, in context of \emph{Practical Statistics}, are:
\begin{enumerate}
	\compresslist
\item{Which prediction methods work best for predicting students' final grades?}
\item{How do we identify a stage in the semester that can adequately balance the required timing of intervention with the quality of the prediction?}
\item{What effect do cluster memberships based on student engagement have on prediction error?}
\end{enumerate}


\section{Method}

In this section we discuss the course background of \emph{Practical Statistics}, as well as the data collection process and analysis used in this study.

\subsection{Course Background Information}
This study took place in University Dublin (UCD). Many of the large first year courses in UCD start in week 1 with material which links to the country's main State Examination and builds from there. Owing to the large class sizes with mixed ability and the progression of material beyond prior knowledge, it may be several weeks before we can identify students who are struggling with the course. \emph{Practical Statistics}, a large online undergraduate course aimed at first years, was selected as an example of a STEM course with weekly continuous assessment. It is designed as an introductory course in statistics for a class of mixed ability students. The lecturer allocates 40\% of the final mark to continuous assessment and distributes the continuous assessment throughout the course semester to encourage students to continuously engage with the course. \emph{Practical Statistics'} lectures are completely online but the students have 24 hours of software labs. The continuous assessment is achieved through: lecture questions based on the course material (weeks 1-12; 0.5\% per week; included in model from week 3); watching all of the online videos (2\%); Minitab lab questions (weeks 3-5; 1\% per week;  included in model at week 5); R lab questions (weeks 7-11 excluding week 8; 1\% per week; included in model at week 11); Minitab lab examination (week 6; 10\%;  included in model at week 6); and R lab examination (week 12; 15\%;  included in model at week 12). Answers to lecture questions and lab sheets are submitted to the LMS and automatically marked by the system, with the marks being returned instantaneously. Students have until midnight of the following Sunday to submit answers.  

\subsection{Participants} 
In the first semester of 2015/16 there were a total of 144 students registered for \emph{Practical Statistics}. Students' data was removed from the study if: students opted out of the research study; students did not take the end of semester examination; or students had personal circumstances which affected how they were officially graded for the course. Students with extenuating circumstances were excluded as these circumstances could impact students' continuous assessment and LMS use. This could impact predictions. In accordance with our ethical permissions from UCD, we removed these students rather than investigating individual student's circumstances. Subsequently our analysis sample included 136 participants from \emph{Practical Statistics}. 

\subsection{Data Collection and Measurements}

Data were recorded for students in regards to three categories: students' background information; continuous assessment; and LMS usage on a fine-grained scale. Background information of students (gender, course type (elective, option or core), registration of students (repeating course etc.), students' year of study, students' programme and Irish/non-Irish) were included as variables to account for differences in educational background and prior experience of students. Online resources (for example videos, lectures slides, pdfs) were grouped into folders based on the material content. In total, there were 15 folders (week 1 course material, ..., week 12 course material, lecture questions solutions, course information, and past examination solutions). For each folder, we included the activity level for the folder for a given week as a variable, for example, in week one student `8979' had an activity level of 12 for the `week 1 course material' (see Table~\ref{Data}). The dataset was designed to be flexible whereby statistical analysis could be performed to incorporate data up to any stage/point in a semester. We performed statistical analysis for the end of each week in the semester (12 teaching weeks) as well as initially (when only background information was available), the end of revision week, and for the end of semester when the written examination was completed. In total this forms fifteen stages. 

\begin{table}[h!]
	\footnotesize
\centering
\caption{Example dataset to be used to predict students' final module mark}
\setlength\tabcolsep{3pt} 
	\begin{tabular}{cccccccc} \hline
		\textbf{\begin{tabular}[c]{@{}c@{}}Student\\ Code\end{tabular}} & \textbf{Gender} & \textbf{...} & \textbf{Major} & \textbf{\begin{tabular}[c]{@{}c@{}}Lecture Q\\ Results\end{tabular}} & \textbf{Week 1 Folder (In Week 1)} & \textbf{Week 1 Folder (In Week 2)} & \textbf{...} \\ \hline
		\textbf{8979} & F & ... & Science & 70.6 & 12 & 23 & ... \\ 
		\textbf{9079} & M & ... & Science & 95.1 & 8 & 15 & ... \\
		\textbf{4567} & M & ... & Arts & 56.8 & 3 & 4 & ... \\
		\textbf{4547} & M & ... & Arts & 64.7 & 7 & 12 & ... \\
		\hline
	\end{tabular}
\label{Data}
\end{table}

\subsubsection{Prediction Methods}

K-fold cross-validation has been used in multiple prediction model studies (Wolff et al., 2013; Azcona \& Casey, 2015). Our prediction models (Random Forest; BART; XGBoost; Principal Components Regression; Support Vector Machine; Neural Network; Multivariate Adaptive Regression Splines; and K-Nearest Neighbours) were run using 10-fold cross-validation for the same folds. The final percentage grade was used as the response variable. 

\begin{itemize}
	\compresslist
\item{Random Forest (RF) is an ensemble learning method. \cite{Breiman2001} state ``random forests are a  combination of tree predictors such that each tree depends on the values of a random vector sampled independently and with the same distribution for all trees in the forest'' (p.~1). For RF regression prediction, the mean prediction of the individual trees is returned.}
\item{\citet{Kapelner2013} explain that BART is a Bayesian approach to nonparametric function estimation using sums of regression trees which allows for flexibility between non-linear interactions. BART differs from other tree ensemble methods (for example RF) owing to the underlying probability model and use of priors. A benefit of this is that we can create confidence intervals for our predicted values.}
\item{XGBoost is a popular scalable machine learning system for tree boosting \citep{Chen2016}. It can handle large datasets as well as sparse matrices. As XGBoost cannot be applied to categorical data, any categorical variables were recoded as binary variables. For XGBoost modelling $\sqrt{n}$ iterations were run where $n$ is the number of variables. XGBoost was applied to 15 stages in the semester (Initially, week 1, ...). Initially XGBoost was used on 18 variables (where categorical variables were transformed to multiple binary variables). The number of variables and iterations increased on a week by week basis as additional Blackboard data became available.}
\item{Principal Components Regression (PCR) is a technique that reduces a high dimensional dataset to a lower dimension dataset and then performing regression. It does this by finding linear transformations of the data whereby the maximal amount of variance is retained \citep{Ilin2010}.}
\item{``Kernel-based learning methods (including Support Vector Machines (SVM)) use an implicit mapping of the input data into a high dimensional feature space defined by a kernel function" \citep*{Karatzoglou2004}. The training of the model is then performed in the feature space.}
\item{Feedforward Neural Network (NN) is a system of nodes which is an imitation of the human brain. A feedforward neural network consists of nodes in layers providing information forward through the layers using the equation $y_{i}=wx_{i}+b$. Training neural networks is considered to be difficult \citep*{Larochelle2009}.}
\item{Multivariate Adaptive Regression Splines is a non-parametric stepwise regression procedure \citep*{Friedman1991} When including variables, the range of the variable is partitioned into subsets and a constant is applied to each subset for regression. In the backward pass, the model is pruned to limit overfitting.}  
\item{K-Nearest Neighbours (KNN) is a nonparametric method whereby the `k' nearest neighbours or `k' most similar cases impact the prediction/classification of the case of interest \citep*{Hechenbichler2004}. In the case of regression, the `k' nearest neighbours response values are averaged with importance weightings being considered}. 
\end{itemize}

We used MAE between the predicted grade and actual grade as a comparison basis to observe the improvement in the accuracy of the prediction model on a week-to-week basis. This allowed us to identify an `optimal time' for an early warning system to be employed. To improve the accuracy of the initial models, our prediction models were applied to different feature sets which included a combination of: continuous assessment data; background information; as well as varying the levels of LMS data. To further reduce the prediction error we considered:  Sunday count variables\footnote{This is discussed further in Section 4.3.}; cumulative count variables; and resultant cluster analysis. The feature sets discussed in Section 4 are:

\begin{enumerate}
	\compresslist
\item{Initial Model - Variables include background information, continuous assessment, and LMS activity level per folder}
\item{No LMS Variables - Variables include background information and continuous assessment}
\item{Cumulative Variables - Variables include background information, continuous assessment, and cumulative activity level for each individual folder (for Sundays and for weekdays)}
\item{Cluster Variables - Variables include background information, continuous assessment, cumulative counter of views for each individual folder (for Sundays and for weekdays), and cluster membership variables}
\end{enumerate}

\subsubsection{Clustering Methods}

The dataset used for clustering contained fine-grained LMS data (the activity level for each individual folder per week and per Sunday). We use the model-based clustering package \emph{mclust} \citep{Scrucca2016} to create an additional clustering of our variables. We use this package because of its repeated superior performance compared to other clustering algorithms \citep{Scrucca2016}, and its ability to model a wide variety of cluster sizes and shapes. Owing to its model-based nature, an advantage of using \emph{mclust} is its ability to calculate probability memberships for each individual to each cluster. Clustering was performed for each stage in the semester. The estimated Bayesian Information Criterion (BIC) was compared for the different combinations, and the combination which maximised the BIC was selected. The resultant cluster membership was considered as a variable for prediction modelling.


\section{Results}

We now describe the development of prediction methods for an early warning system. Continuous assessment played an important role in our modelling. When developing an early warning system, we need to account for any delays in the correction of continuous assessment or collection of data for example if a midterm in week 5 takes two weeks to correct, we should include it in week 7. \emph{Practical Statistics} benefits from the instantaneous nature of online LMS assignments. Through the development of our early warning system, we are able to identify an optimal time (week 5-6) in the \emph{Practical Statistics'} semester to apply an early warning system. 

\subsection{Student Engagement and Continuous Assessment}

\begin{figure}[ht!]
	\centering
	\graphicspath{C:/Users/SMSuser/Documents/}
	\includegraphics[width=160mm,height=80mm]{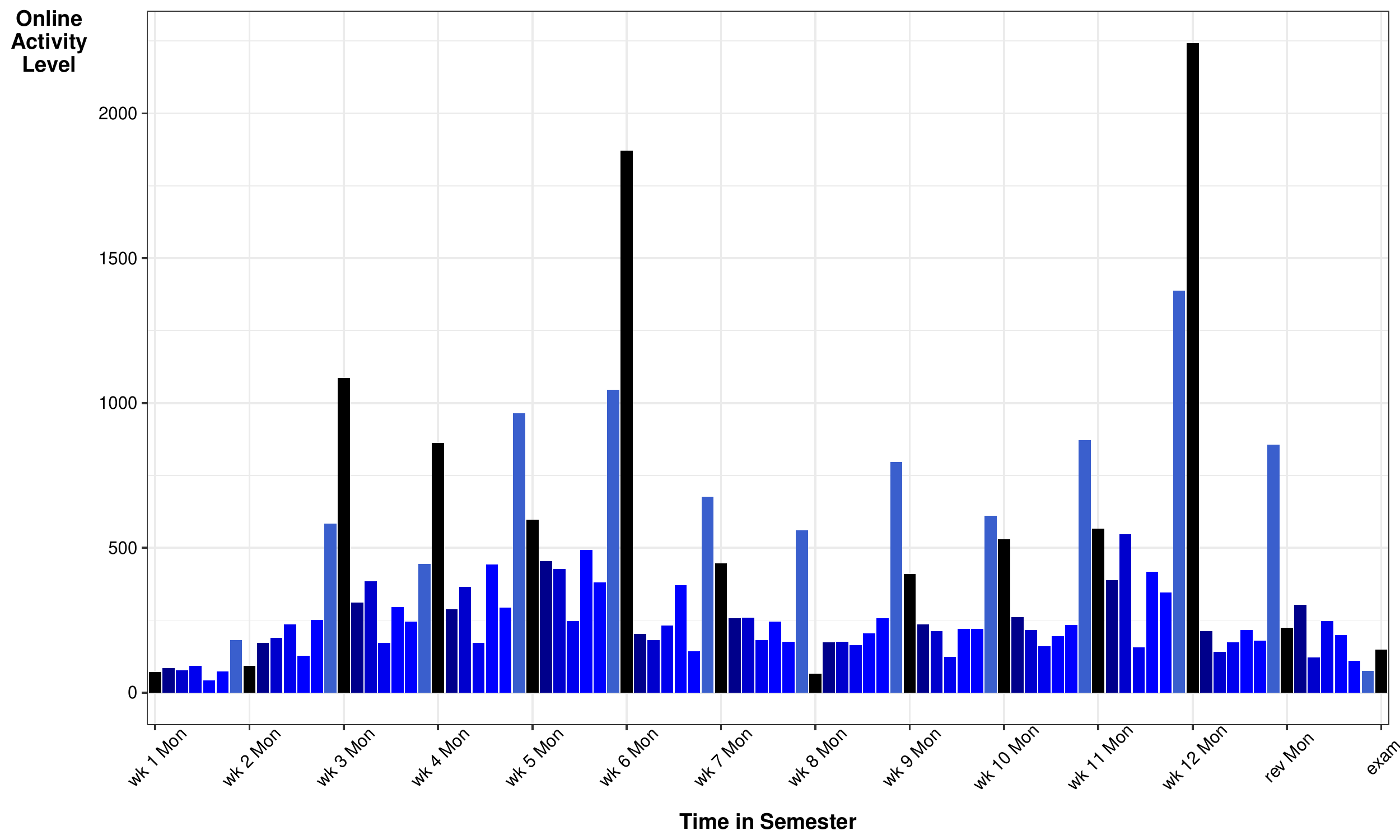}
	\caption{Activity level of online resources per day over the semester}
	\label{Engagement}
\end{figure}

\citet*{Holmes2015} and \citet*{Cole2012} have suggested that continuous assessment encourages student engagement. As previously mentioned, \emph{Practical Statistics} was designed to ensure consistent student engagement through having continuous assessment on a weekly basis throughout the course. Figure~\ref{Engagement} shows that online materials were accessed throughout the semester, however the level of activity, not surprisingly, varied across the semester.

The deadline for weekly online lecture questions for credit was on Sunday nights, and this corresponds with the weekly peak in online resource activity. These peaks might suggest two types of students: students who study immediately prior to assessments; and students who study in advance of assessment. Similarly, as expected, the time when the greatest number of resources was accessed corresponds to the day of the R lab examination (the Monday of week 12). A similar peak occurs on the Monday of the Minitab lab examination in week 6. This connection between online views and continuous assessment suggests that a key driver of students' interaction with online resources is continuous assessment. 

\subsection{Clustering Analysis}

\emph{mclust} was applied to several variations of the dataset. Considering the high number of view counts on Sunday, this included investigating the potential of Sunday online activity as separate to weekday\footnote{In this study weekday view counts includes Saturday view counts.} online activity. After investigating resultant clusters, \emph{mclust} was only applied to fine-grained LMS data (the activity level for each individual folder per week and per Sunday). Continuous assessment variables and background information of students, were not included as cluster variables. The resultant clusters identified differences in students' frequency levels of using online resources. In comparison to \citet{Lust2011} who divides online resources into tool types, this method is cruder as the clustering is unlikely to pick up subtle differences in students' learning strategies.

\begin{figure}[h]
	\centering
	\footnotesize
	\graphicspath{C:/Users/SMSuser/Documents/}
	\includegraphics[width=138mm,height=93mm]{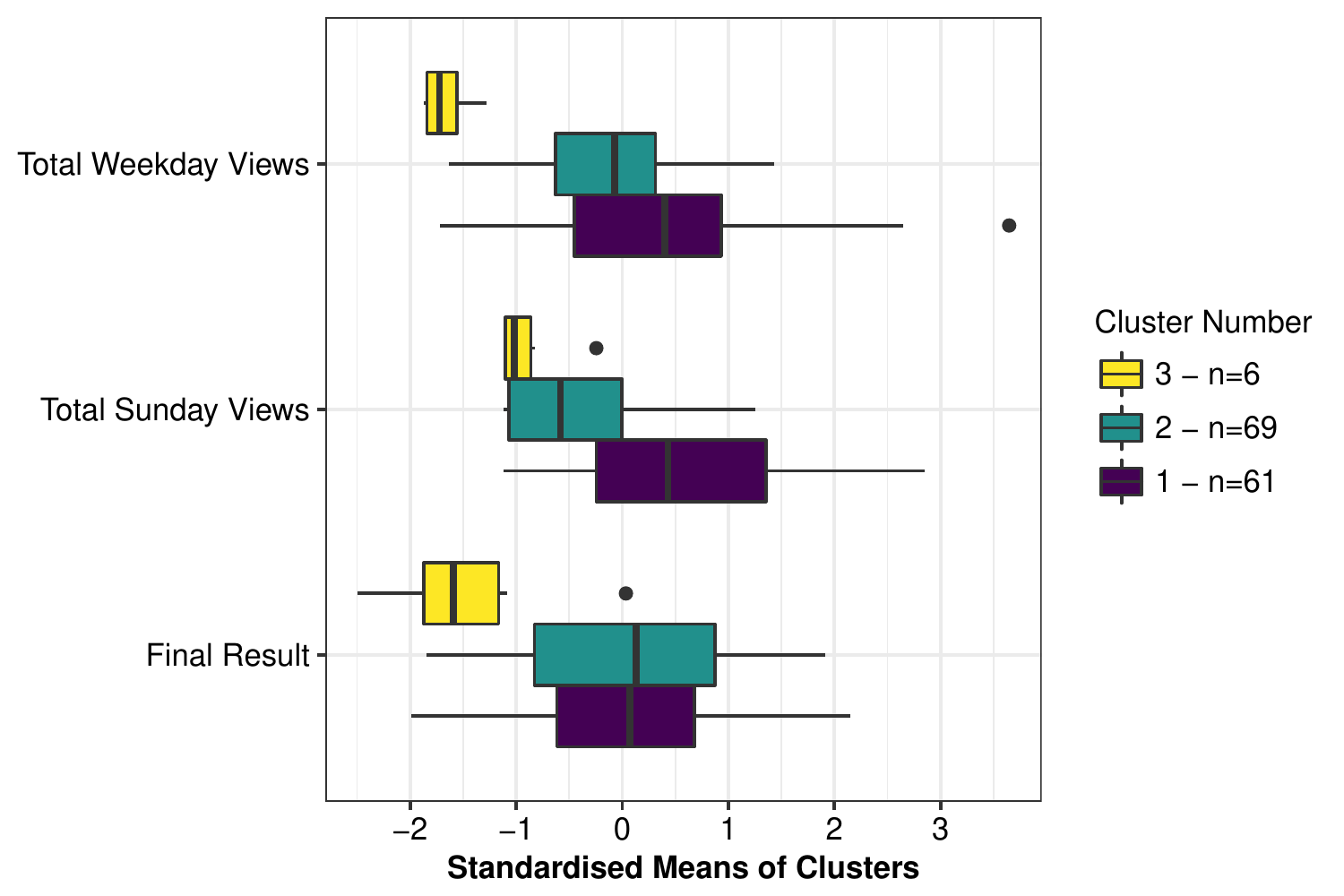}
	\caption{Identifying engagement patterns of Practical Statistics through boxplots of selected standardised variables for week 5. for example cluster 3 contains six students who have below average resource usage.}
	\label{Identifying Engagement Patterns of Practical Statistics}
\end{figure}

For example, for week 5 (identified optimal time) the variables used were fine-grained LMS data (the activity level for each folder per week and per Sunday) for weeks 1-5. \emph{mclust} identified 3 clusters ($n_1=61 $ $(44.9\%)$, $n_2=69$ $(50.7\%)$ and $n_3=6$ $(4.4\%)$). The distinct clusters are best represented in 2D format by boxplots (see Figure~\ref{Identifying Engagement Patterns of Practical Statistics}) showing the standardised means and spread of the selected variables for each cluster. Three variables (Total Weekday Views (up to week 5), Total Sunday Views (up to week 5), and Final Grade) were selected to show the distinct clusters (see Figure~\ref{Identifying Engagement Patterns of Practical Statistics}). 
 
Cluster 3 students are students who display below average engagement with resources and have the widest final grade range. \citet{Lust2011} categorize these as no-users or low frequency users. Cluster 2 represents the students who have below average resource use on Sunday, and average resource use during the week. In comparison, Cluster 1 represents students who engage above average with resources overall. Despite this high engagement, they have the median final grade. \citet{Lust2011} would describe these as the intensive users. Subsequently, as the cluster analysis displayed distinct clusters with various engagement patterns, students' cluster group membership was used as variables in the prediction analysis.

\subsection{Prediction Modelling}

\begin{figure}[h]
	\centering
	\graphicspath{C:/Users/SMSuser/Documents/Developing Early Warning Systems/}
	\includegraphics[width=132mm, height=104mm]{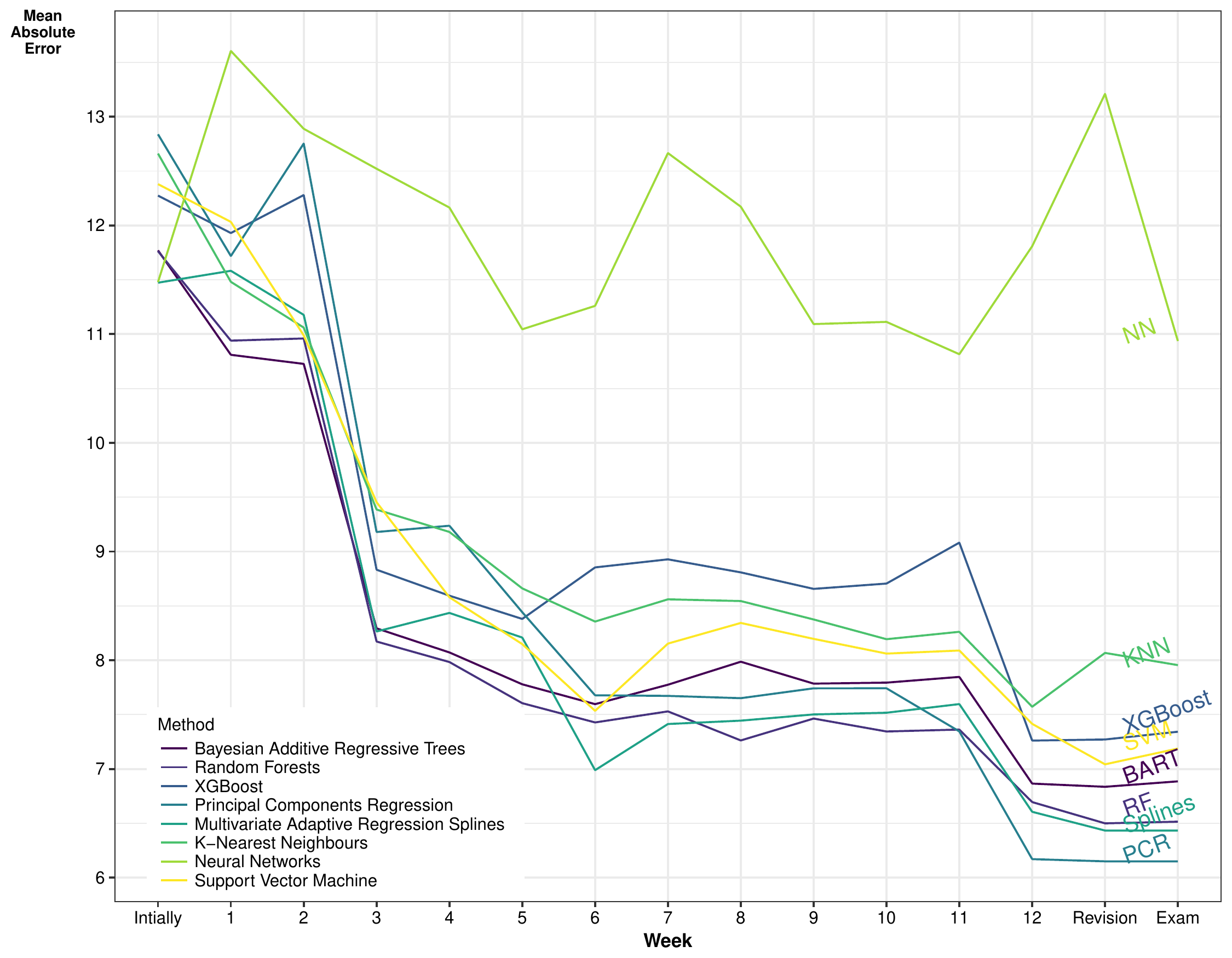}
	\caption{Average MAE per student on a week-by-week basis from multiple out of sample prediction methods}
	\label{Initial Model}
\end{figure}

Initial prediction modelling was performed on the dataset for each week (all variables available up to that date were included - see Initial Model Section 3.3.1) to determine an optimal time for corrective measures. The initial stage (before teaching semester began) and final stage of the semester acted as a baseline for comparison for the power of the prediction model (see Figure~\ref{Initial Model}). Out of the methods investigated, Neural Networks is clearly the inferior method.

In Figure~\ref{Initial Model} Principal Components Regression achieves the lowest MAE value  (approximately 6 points at week 12). PCR reduces the number of variables before performing regression. An interesting feature of Figure~\ref{Initial Model} is the substantial decrease in error from week 2 to week 3. This decrease in error coincides with the inclusion of continuous assessment in the prediction model (the deadline for week 1 lecture questions was in week 3). This emphasises the role continuous assessment plays as a predictor in online STEM courses. To confirm the importance of continuous assessment we investigated the variable importance of the models. Every model selected continuous assessment variables as the main variables in the model. Figure~\ref{Initial Model} also shows that between weeks 7 and 11 there is relatively little change in the predictive power of the models. Subsequently, the stages up to week 6 were identified as important for further data analysis. For early warning systems, a balance is required between the accuracy of the prediction models and the stage in the semester. The stage in the semester needs to reflect where corrective measures could most effectively be given to students. For \emph{Practical Statistics}, week 5 is potentially the optimal time for implementing an early warning system. We have included our R code for this in the Appendix with more detailed R code and fictitious datasets available on GitHub at https://github.com/ehoward1/Early-Warning-System. 

\begin{figure}[h]
\centering
\graphicspath{C:/Users/SMSuser/Documents/}
\includegraphics[width=140mm,height=110mm]{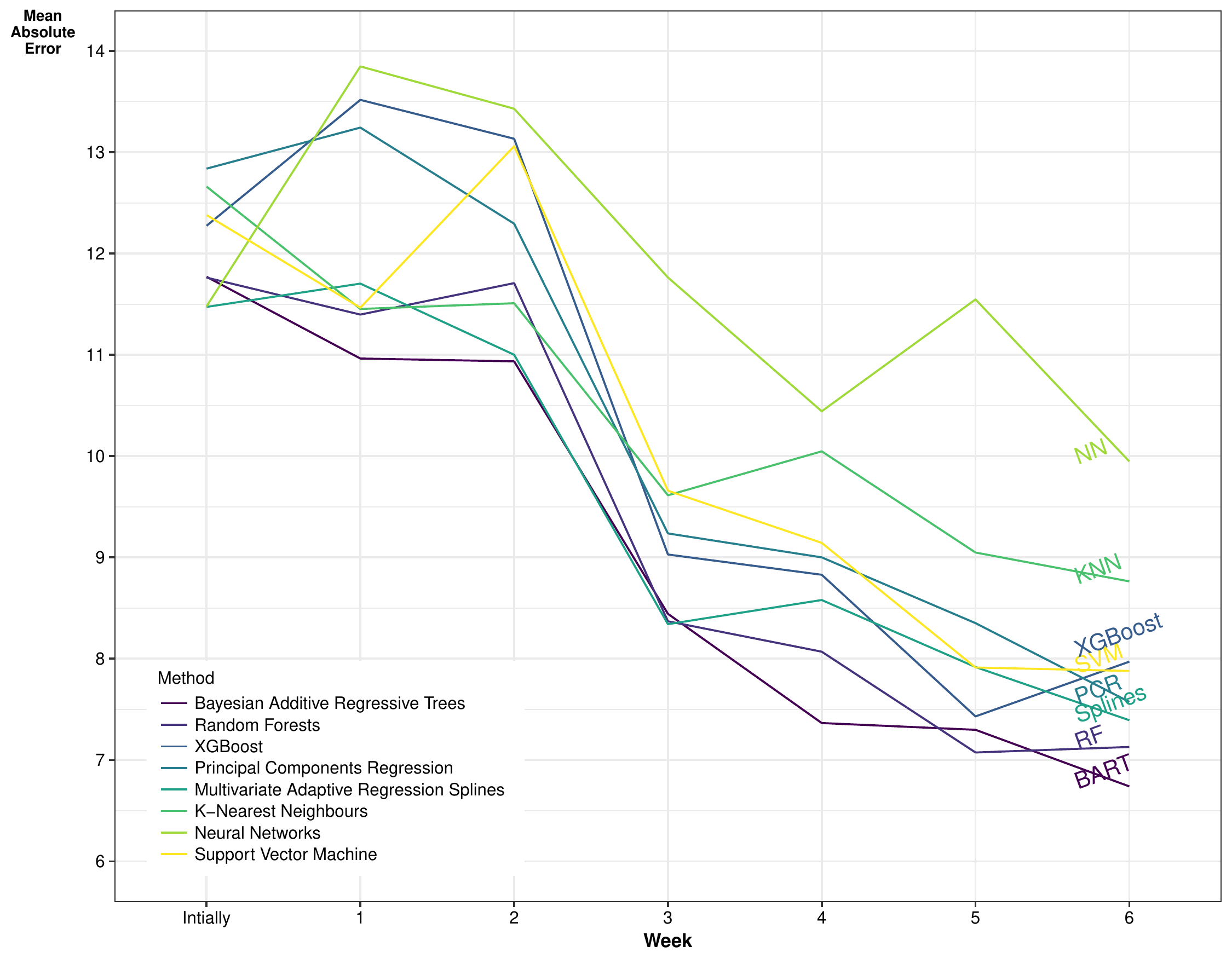}
\caption{Improving the Prediction Model by changing the feature dataset to include cumulative variables for LMS interactions and cluster membership variables}
\label{Improving Models}
\end{figure}

We considered alternative feature sets including removal of the LMS data (which provided slightly less accurate predictions), including cumulative activity level for each folder (Cumulative Variables dataset), and including cluster membership variables. Progressing, we will look at the Cluster Variable dataset in further detail. The Cluster Variables dataset for each student consists of: background information; continuous assessment; and cumulative counter for the activity level of each individual folder (for Sundays and for weekdays) as well as cluster membership variables. While including the cluster variable (in most cases) does not alter the MAE significantly, clustering can provide us with information about student engagement in general which may be of value (see Section 4.2 Clustering Analysis). Figure~\ref{Improving Models} gives the average MAE per student for the Cluster Variables dataset up to the optimal time of week 6. The second substantial decrease in MAE between weeks 4 and 5 corresponds to the second inclusion of continuous assessment (Minitab lab results). Using our BART predictive model we can identify the final mark the student will obtain to approximately a MAE of 6.5 at week 6. We will proceed by discussing in further detail the BART prediction model at week 5 using the Cluster Variable dataset. This dataset consists of 29 explanatory variables.

\begin{figure}[ht!]
\centering
\graphicspath{C:/Users/SMSuser/Documents/}
\includegraphics[width=150mm,height=100mm]{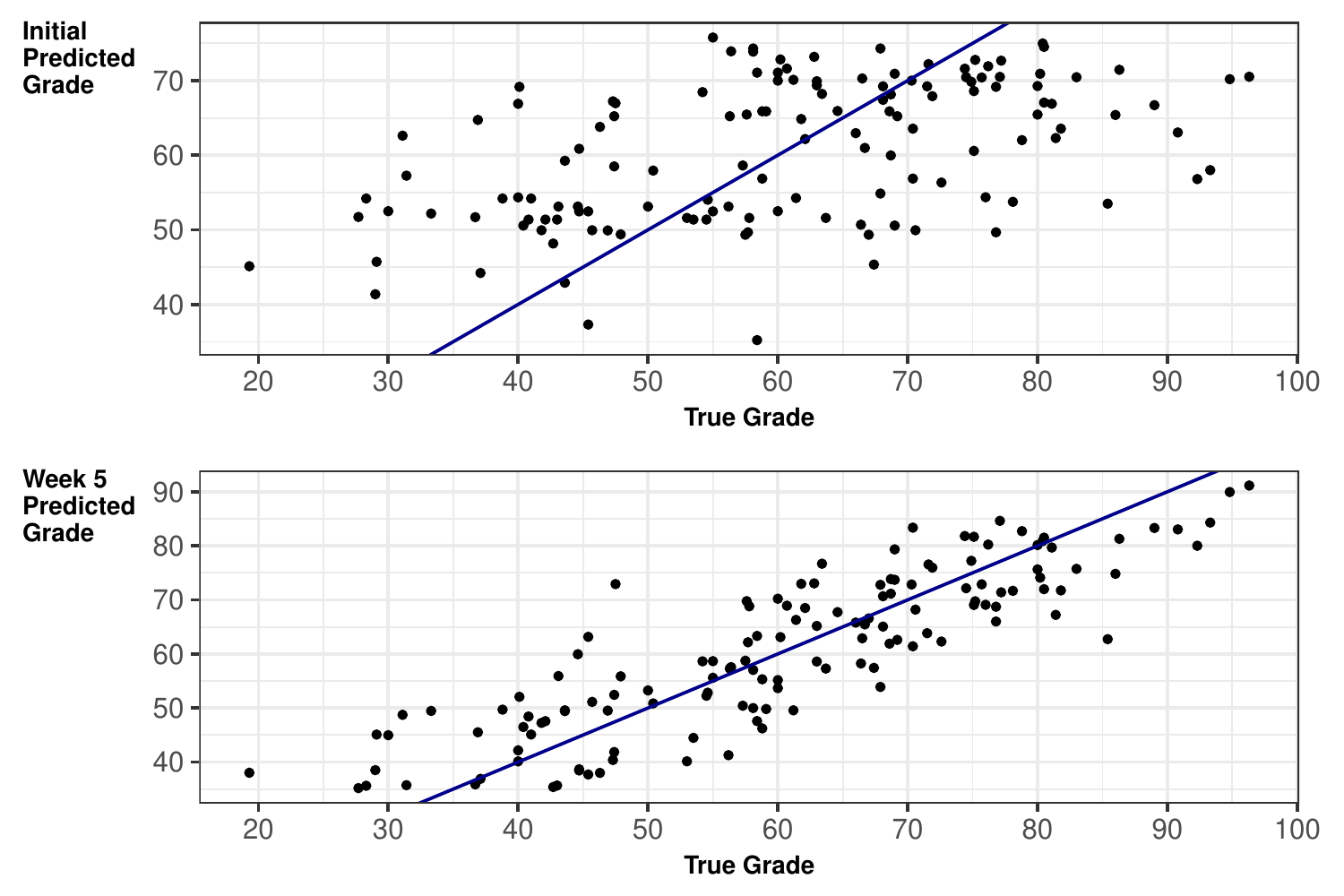}
\caption{Scatter Plots showing Predicted Grade versus Actual Grade Initially and at Week 5 via an Out of Sample 10-Fold Cross-Validation}
\label{Scatter Plots}
\end{figure}

We can visually determine the performance of our predictive model by plotting the predicted final grade against the true final grade for each student. Figure~\ref{Scatter Plots} shows the predicted grade plotted against the actual grade of each student, both initially and at the end of week 5. An identity line, showing when the predicted grade equals the actual grade (i.e. a perfect prediction), has been included in Figure~\ref{Scatter Plots}. The initial plot acts as a baseline, displaying how the initial prediction of final grade has very low correlation with the actual grade of students i.e. a poor predictive performance. The initial model relies on a limited number of background/demographic variables. As several students have the same background information, this has resulted in multiple students receiving the same predicted grade. This has resulted in `bands' of predicted grades. In comparison, the second plot's data is quite linear ($R^2$ = 0.74) and tighter to the identity line,  with some outliers. It suggests that by week 5 we can make reasonable grade predictions as the grade predictions are strongly correlated to the actual grade. This supports the belief that week 5 is an optimal time to implement an early warning system, and that the selected BART model (Method - Cluster Variables) performs competently. 

\section{Discussion and Conclusion}

\subsection{Continuous Assessment}
 
The variables used in this study were divided into three categories: students' background information; students' engagement with LMS; and continuous assessment results. Continuous assessment proved unsurprisingly the most important category. Continuous assessment variables were repeatedly chosen as the most important variables by all of the prediction models. Continuous assessment encourages students to engage with a course (Holmes, 2015) and partially accounts for the different levels of LMS interaction throughout the semester. This is observable from the spikes in LMS resource use prior to continuous assessment tests and deadlines (Figure~\ref{Engagement}).  We suggest the inclusion of consistent continuous assessment in online courses encourages students' engagement over the entire semester (as stated by \citet*{Cole2012}), and limits the number of students studying only in the weeks prior to the final examination. The addition of continuous assessment also contributes to minimising prediction error when building early warning systems however this should not be the main reason for its inclusion. We hypothesize that a low percentage for continuous assessment would also achieve the same effect provided that the continuous assessment is throughout the semester.

This study investigates how to approach developing an accurate prediction model for an early warning systems. The dataset which only had continuous assessment and background information variables, performed comparatively well to the other feature sets, and enjoyed the benefit of being the simplest model. However, by using this model we fail to identify areas of the curriculum where students struggled. A key element in learning analytics is using the resultant analysis for the benefit of the student and teacher. By including the extra LMS variables, we are able to investigate for individual students aspects of the curriculum that they failed to engage with or had overly high engagement with (potentially a sign of a harder concept or an area with which the student struggled). This advantage for the inclusion of LMS variables is considerable, and should be weighed against the simplicity of the `No LMS Variables' data set.

\subsection{Advancements in Developing Early Warning Systems}

This study summarises the methods employed in developing prediction models for early warning systems, and builds upon the current work. Unlike other studies, we do not dichotomise students' final marks to pass/fail. We discuss how one may reduce the prediction error through: use of fine-grained variables; manipulation of variables; and the inclusion of cluster membership in prediction modelling. The detail provided by fine-grained variables gives more information on students' engagement patterns. Subsequently, we hypothesize that analysis of fine-grained variables will allow for more personalised corrective measures. We have used predictive methods (BART and XGBoost) which are uncommon in the data analytics literature as well as common predictive methods (Neural Networks, K-Nearest Neighbours and Random Forest). We found that decision tree methods perform particularly well (BART and Random Forest). Decision tree methods are suitable when using a large number of variables. Hence BART, a decision tree method, is appropriate when using fine-grained variables. BART may be preferable over other decision tree methods, for example Random Forests, owing to its Bayesian nature which allows for the inclusion of error variance which is independent of tree structure and leaf parameters \citep{Kapelner2013}. In our study BART outperformed the other prediction models tested at the optimal time of weeks 5-6.

Clustering is not a necessary step in developing prediction models. However, we have shown that clustering can be used to identify distinct student patterns of engagement which can be used to further reduce the prediction error. Also, clustering may help to identify how students approach learning and subsequently be used to provide corrective measures. The method outlined in this study is appropriate for both online courses and large classes with a significant amount of online material. Through combining these methods, we obtain an average prediction error (based on out of sample 10-fold cross validation and MAE) of 6.5 percentage points by week 6.  

A key part of this study was identifying an `optimal time' to implement an early warning system. Implementing an early warning system too early would result in inaccurate identification of (at-risk) students. In contrast, implementing it too late would diminish the effect of supporting and helping students. Data analysis of prediction models identify week 5/6 as the critical time in the semester for \emph{Practical Statistics} whereby prediction models have reasonably accurate forecasts balanced with sufficient time to intervene and support at-risk students. Identifying at-risk students is only one stage in an early warning system, another stage is understanding what effective supports should be provided to students.  Consequently, our current research involves identifying at-risk students during the `optimal time' in \emph{Practical Statistics} and examining which feedback/intervention measures are effective for large STEM courses.  

\subsection{Limitations}

The method outlined in this study discusses how to develop an accurate predication model for an early warning system for a course, and how to recognise an optimal time to provide students with corrective measures during a course. \emph{Practical Statistics} is an example of a STEM course which has continuous assessment distributed weekly throughout the semester. The method discussed in this study may not be an optimal method for other online courses, particularly if the course is from a significantly different academic field. Each course is unique and will have its own unique feature set. STEM based courses, particularly early undergraduate courses, tend to have continuous assessment which ties to the final examination. We believe BART is applicable for these STEM courses. 

For the purpose of reproducibility, the R code for comparison of the prediction models has been included in Appendix A. Further code for for this study is available on GitHub at https://github.com/ehoward1/Early-Warning-System with fictitious datasets (owing to ethical constraints). 

\section*{Acknowledgements}
We would like to thank UCD IT services for providing us with Blackboard data. This research did not receive any specific grant from funding agencies in the public, commercial, or not-for-profit sectors.

\section*{Ethics}
This study was conducted in accordance with UCD ethics guidelines, and approved by the UCD Ethics Committee under application number LS-15-53-Meehan. 

\clearpage

\section*{Appendix - R code}

Function to run and compare 10-cross fold validation for all prediction methods used in the paper. Fictitious datasets for this function and further R code for this paper are available on GitHub at https://github.com/ehoward1/Early-Warning-System. 

\begin{verbatim}
require(xgboost)
require(randomForest)
require(bartMachine)
require(pls)
require(caret)
require(magrittr)
require(earth)
require(nnet)
require(car)
require(kknn)
require(kernlab)

prediction_function <- function(dataset, dataset_boost){ # dataset boost is for xgboost 

	set.seed(123)
	folds = createFolds(1:nrow(dataset), k = 10, list = FALSE)
	dataset_boost = apply(dataset_boost, 2, as.numeric) # XGBoost runs for numeric data, not integers

	# Vectors to store error for each prediction methods
	pred_bm = vector("numeric")
	pred_rf = vector("numeric")
	pred_pcr = vector("numeric")
	pred_xg = vector("numeric")
	pred_kknn = vector("numeric")
	pred_svm = vector("numeric")
	pred_nnet = vector("numeric")
	pred_earth = vector("numeric")
	grades = vector("numeric")

	# Loop through the folds for each prediction method
	for(i in 1:10){

		# Setting up the data
		train = dataset[folds!=i,] %>% data.frame %>% na.omit
		train_b = dataset_boost[folds!=i,] %>% data.frame %>% na.omit 
		test = dataset[folds==i,] %>% data.frame
		test_b = dataset_boost[folds==i,] %>% data.frame 

		# BART
		bm = bartMachine(train[,-1], train[,1], seed = 123, alpha = 0.95, num_burn_in = 400,
		num_tree = 100, num_rand_samps_in_library = 20000, k = 2, q = 0.9, nu = 3)
		pred_bm = c(pred_bm, predict(bm, test[,-1]))

		# Random Forest (RF)
		rf = randomForest(train[,-1], train[,1], ntree = 100)
		pred_rf = c(pred_rf, predict(rf, test[,-1]))

		# Principle Components Regression (PCR)
		pcr = pcr(FINAL~., data = train)
		var_exp = compnames(pcr, explvar = TRUE)
		var_e = unlist(strsplit(var_exp, "[ (]")) %>% as.numeric() %>%  setdiff(c(1:150, NA))
		var_total = 0

# Calculating number of variables to include based on variation explained
for(j in 1:length(var_e))
{
var_total = var_total + var_e[j]
if(var_e[j] < 1 || var_total > 90)
{
n_comp = j
break
}
}
pred_pcr = c(pred_pcr, predict(pcr, test[,-1], ncomp = n_comp))

# Xgboost
iter = train_b %>% ncol %>% sqrt %>% ceiling
xg = xgboost(data = as.matrix(train_b[,-1]), label = train_b[,1], eta = 0.5, 
nround = iter, max.depth = 4, objective = "reg:linear")
pred_xg = c(pred_xg, predict(xg, as.matrix(test_b[,-1])))

# K-Nearest Neighbours (KNN)
kknn = train.kknn(FINAL ~., kmax = 15, distance = 1, data = train)
pred_kknn = c(pred_kknn, predict(kknn, test[,-1]))

# Neural Network (NN)
my.grid = expand.grid(.decay = c(0.05, 0.5, 0.75), .size = c(4, 9))
nnet = train(FINAL~., data = train, linout = 1, 
method = "nnet", maxit = 500, tuneGrid = my.grid, trace = FALSE) 
pred_nnet = c(pred_nnet, predict(nnet, test[,-1]))

# Support Vector Machine (SVM)
svm = ksvm(FINAL ~., data = train, C = 5)
pred_svm = c(pred_svm, predict(svm, test[,-1]))

# Multivariate Adaptive Regression Splines (Splines)
earth = train(FINAL~., data = train, method = "earth",
tuneGrid = data.frame(degree = c(1,2), nprune = 5)) 
pred_earth = c(pred_earth, predict(earth, test[,-1]))

grades = c(grades, test$FINAL)
}

# Calculating the error for each method
error_rf = sum(abs(pred_rf - grades))/nrow(dataset) 
error_pcr = sum(abs(pred_pcr - grades))/nrow(dataset) 
error_xg = sum(abs(pred_xg - grades))/nrow(dataset_boost) 
error_bm = sum(abs(pred_bm - grades))/nrow(dataset) 
error_earth = sum(abs(pred_earth - grades))/nrow(dataset) 
error_kknn = sum(abs(pred_kknn - grades))/nrow(dataset) 
error_nnet = sum(abs(pred_nnet - grades))/nrow(dataset) 
error_svm = sum(abs(pred_svm - grades))/nrow(dataset) 

# Returning Values
my_list = list("MAE_bm" = error_bm, "MAE_rf" = error_rf, "MAE_pcr" = error_pcr, 
"MAE_xg" = error_xg, "MAE_kknn" = error_kknn, "MAE_nnet" = error_nnet, 
"MAE_svm" = error_svm, "MAE_earth" = error_earth)   

return(my_list)
}

\end{verbatim}

\clearpage

\bibliography{Bibliography}
\bibliographystyle{chicago}

\end{document}